\documentclass[oneside,a4paper]{article}
\usepackage{globalization}

\begin{document}

\title{A globalization for non-complete\\ 
but geodesic spaces}
\author{Anton Petrunin}
\date{}
\maketitle

\begin{abstract}
I show that if a geodesic space has curvature bounded below locally in the sense of Alexandrov then its completion has the same lower curvature bound globally.
\end{abstract}

\section*{Introduction}

Let us recall few definitions.

\begin{itemize}
\item A \emph{geodesic space} is a metric space such that any two points can be connected by a minimizing geodesic.
\item A \emph{length space} is a metric space such that any two points can be connected by a curve with length arbitrary close to the distance between the points.
\item A complete length space is an 
\emph{Alexandrov space with curvature $\ge \kappa$} 
if for any quadruple $(p;x^1,x^2,x^3)$ the \emph{(1+3)-point comparison} holds;
i.e.,  
if 
\[\angk\kappa p{x^1}{x^2}
+\angk\kappa p{x^2}{x^3}
+\angk\kappa p{x^3}{x^1}\le 2\cdot\pi\]
or at least one of the model angles $\angk\kappa p{x^i}{x^j}$ is not defined.
\end{itemize}

\parbf{Theorem.}
\textit{Let $X$ be a geodesic space.
Assume that for any point $x\in X$
there is a neighborhood $\Omega\ni x$ 
such that the $\kappa$-comparison holds for any quadruple of points in $\Omega$.
Then the completion of $X$ is an Alexandrov space with curvature $\ge \kappa$.}
\medskip

The question was asked 
by Victor Schroeder around 2009.
Later I learned from Stephanie Alexander,
that this statement has an application.

I want to thank 
Stephanie Alexander, 
Richard Bishop,
Urs Lang,
Nan Li,
Victor Schroeder  
for the interesting discussions.
A special thanks to Alexander Lytchak 
for correcting few mistakes
and suggesting an easier proof of Lemma~\ref{lem:xyz} via his favorite ultrapower.

\medskip

In Section~\ref{sec:start} I give a reformulation of the above theorem, which will be proved in Section~\ref{sec:end} 
after auxiliary statements proved in Section \ref{sec:prelim}.

\section{Observation and reformulation.}
\label{sec:start}

The distance between points $x$ and $y$ in a metric space will be denoted as 
$\dist{x}{y}{}$.

A geodesic from $x$ to $y$ will be denoted as $[xy]$,
once we write $[xy]$ we implicitly make a choice of one such geodesic.
Also we use the following short-cut notation:
\begin{align*}
\l] x y \r]&=[xy]\backslash\{x\},
&
\l] x y \r[&=[xy]\backslash\{x,y\},
&
\l[ xy \r[&=[xy]\backslash\{y\}.
\end{align*}
We denote by $\MM^\kappa$ the model $\kappa$-plane
and 
$\varpi^\kappa\df\diam\MM^\kappa$ (so $\varpi^\kappa=\infty$ if $\kappa\le 0$ and $\varpi^\kappa=\tfrac\pi{\sqrt{\kappa}}$ otherwise).

\parbf{Angles.}
Assume $[px]$ and $[py]$ be two geodesics in $X$
and $\bar x\in \l]px\r]$ and $\bar y\in \l]py\r]$.
Since $\kappa$-comparison holds in a neighborhood of $p$
the function
$$(\dist{p}{\bar x}{},\dist{p}{\bar y}{})
\mapsto
\angk\kappa p{\bar x}{\bar y}$$
is nonincreasing in both variables 
for sufficiently small values of the variables.
It follows that angle 
$$\angle\hinge pxy
\df
\lim
\set{\angk\kappa  p{\bar x}{\bar y}}
{\bar x\in \l]px\r],\bar y\in \l]py\r],\dist{p}{\bar x}{}\to 0,\dist{p}{\bar y}{}\to 0}$$
is defined for any hinge $\hinge pxy=([px],[py])$ in $X$.

The same way as for Alexandrov space,
one can show that 
$$\angle\hinge pxy+\angle\hinge pyz+\angle\hinge pzx\le 2\cdot\pi$$
for any three hinges formed by geodesics $[px]$, $[py]$ and $[pz]$.
It follows that, if $p\in \l]xy\r[$ then
$$\angle\hinge pyz+\angle\hinge pzx=\pi.\eqlbl{bigstar}$$


\parbf{Reformulation.}
To prove the theorem, it is sufficient to show that 
for any $\kappa_1<\kappa$, any point $p$ and any geodesics $[qs]$ in $X$ we have 
$$\angk{\kappa_1} q{\bar s}p
\le
\angle\hinge q{\bar s}p,
\eqlbl{({*})}$$
if 
$\bar s\in \l]qs\r]$
and $\dist{q}{\bar s}{}>0$ is small enough.

Indeed, once it is proved, it follows that the inequality \ref{({*})} holds for all $\bar s\in \l]qs\r]$.
Together with \ref{bigstar},
it implies that the (1+3)-point comparison for all quadruples in $X$;
this can be done exactly the same way as in Alexandrov space, see \cite[2.8.2]{bgp}.

Hence the completion $\bar X$ of $X$
is an Alexandrov space with curvature $\ge\kappa_1$ for any $\kappa_1<\kappa$.
From the standard globalization theorem (as it is stated in \cite{akp}) we get the result.

\section{Auxiliary statements}\label{sec:prelim}

As above we denote by $\bar X$ the completion of $X$.

Note that any point $x\in X$ admits an open neighborhood $\Omega$ in $\bar X$ such that $\kappa$-comparison holds for any quadruple of points in $\Omega$.
In particular the following condition holds for any point $p\in\Omega\cap X$ and geodesic $[qs]$ which lie in $\Omega\cap X$.
$$\angk{\kappa_1} q{\bar s}p
\le
\angle\hinge q{\bar s}p,\eqlbl{(**)}$$
where 
$\bar s\in \l]qs\r]$.

An open domain $\Omega$ in  $\bar X$ 
which satisfy \ref{(**)}
will be called a \emph{$\kappa$-domain}.

Note that to prove that $\Omega$ is a $\kappa$-domain,
it is sufficient to check that \ref{(**)}
holds only if $\dist{q}{\bar s}{}$ is sufficiently small.
I.e., if for any point $p$ and geodesic $[qs]$ in $\Omega\cap X$, the condition \ref{(**)} holds if $\dist{q}{\bar s}{}$ is small enough then \ref{(**)} holds for all $\bar s\in\left]qs\right]$.


Note also that if $B(p,2\cdot R)$ is a $\kappa$-domain in $\bar X$ then 
$\kappa$-comparison holds for any quadruple of points in $B(p,R)$.
The later is proved exactly the same way as in Alexandrov space:
for a quadruple $(p;x,y,z)$ we choose a geodesic $[px]$ and apply \ref{bigstar} together with \ref{(**)}  for $\bar x\in \l]px\r[$ such that $\bar x\to p$.
(Everything works since geodesics with ends in $B(p,R)$ can not leave $B(p,2\cdot R)$.)

In particular, if $\Omega$ is a $\kappa$-domain in $\bar X$
then the curvature at each point of $\Omega$ is $\ge \kappa$.
Therefore any local construction in Alexandrov geometry 
can be performed inside of $\Omega$.

For example, we can construct so called \emph{radial curves} 
as far as we do not get out of $\Omega$.
The radial curves are formed by trajectories which try to escape from a given point $w$ using the greedy algorithm; 
these curves parametrized in a special way which makes them behave as geodesics in terms of comparisons.

The following proposition is a local version of \emph{radial monotonicity} in \cite{akp}
and can be proved exactly the same way.

\begin{thm}{Proposition}\label{prop:rad-loc}
Let  
$\Omega\subset \bar X$ be a $\kappa$-domain
and $w,a\in \Omega$.
Assume that $\bar B[w,R]\subset \Omega$ 
and 
$$\dist{a}{w}{}=r<R<\tfrac{\varpi^\kappa}2.$$
Then there is a radial curve $\alpha\:[r,R]\to \Omega$ with respect to $w$
such that $\alpha(r)=a$
and the distance $\dist{p}{\alpha(t)}{}$ 
satisfies the radial monotonicity 
for any point $p$ in $\Omega$.

I.e., if $[\~w\,\~p\,\~\alpha(t)]$ is 
a triangle in $\mathbb M[\kappa]$ 
with sides
\begin{align*}
\dist{\~w}{\~p}{}&=\dist{w}{p}{}
&
\dist{\~p}{\~\alpha(t)}{}&=\dist{p}{\alpha(t)}{}
&
\dist{\~w}{\~\alpha(t)}{}&=t
\end{align*}
then the function 
$$t\mapsto \angle\hinge{\~w}{\~v}{\~\alpha(t)}$$ 
is a nonincreasing in its domain of definition.
\end{thm}

The proposition above is used to prove the following lemma.

\begin{thm}{Key Lemma}\label{lem:key}
Let $\Omega_p$ and $\Omega_q$ be two $\kappa$-domains in $\bar X$.
Let 
\begin{align*}
p&\in X\cap \Omega_p,
& 
q&\in X\cap\Omega_q,
&
w&\in X\cap\Omega_p\cap\Omega_q.\end{align*}

Consider a triangle
$[\~p\~w\~q]$ in $\MM^\kappa$ 
such that 
\begin{align*}
\dist{\~p}{\~w}{}&=\dist{p}{w}{}
&
\dist{\~q}{\~w}{}&=\dist{q}{w}{}
&
\angle\hinge{\~w}{\~p}{\~q}&=\angle\hinge{w}{p}{q}
\end{align*}
Set $R$ to be the distance from the side $[\~p\~q]$
to $\~w$.

Assume $R<\tfrac{\varpi^\kappa}2$
and $\bar B[w,R]\subset \Omega_p\cap\Omega_q$.
Then 
$$\dist{p}{q}{}\le \dist{\~p}{\~q}{}.$$

\end{thm}

\parit{Proof of the Key Lemma.}
Note that if $R=0$, then the lemma follows from the triangle inequality;
further we assume $R>0$.

Let $\~a\in [\~p\~q]$  be a point
which minimize the distance to $\~w$;
so $R=\dist{\~w}{\~a}{}$.

Fix small $\delta>0$;
any value $\delta<\tfrac1{10}\cdot\min\{1,\~R/{\dist{w}{p}{}},{\~R}/{\dist{w}{q}{}}\}$
will do.
Choose $p_\delta\in\l]wp\r]$
and $q_\delta\in\l]wq\r]$ 
so that 
\begin{align*}
\dist{w}{p_\delta}{}&=\delta\cdot\dist{w}{p}{},
&
\dist{w}{q_\delta}{}&=\delta\cdot\dist{w}{q}{}.
\end{align*}
Note that geodesic 
$[p_\delta q_\delta]$ lies in $\bar B[w,\~R]$.
By Alexandrov's lemma (stated as in \cite{akp}), 
one can choose a point $a_\delta\z\in [p_\delta q_\delta]$
so that
\begin{align*}
\angk\kappa w{p_\delta}{a_\delta}
&\le \angle\hinge {\~w}{\~p}{\~a},
&
\angk\kappa w{q_\delta}{a_\delta}
&\le \angle\hinge {\~w}{\~q}{\~a}.
\intertext{Note that} 
\angk\kappa w{p}{a_\delta}&\le\angk\kappa w{p_\delta}{a_\delta},
&
\angk\kappa w{q}{a_\delta}&\le\angk\kappa w{q_\delta}{a_\delta}.
\intertext{Therefore}
\angk\kappa w{p}{a_\delta}
&\le \angle\hinge {\~w}{\~p}{\~a},
&
\angk\kappa w{q}{a_\delta}
&\le \angle\hinge {\~w}{\~q}{\~a}.
\end{align*}

\begin{wrapfigure}{r}{50mm}
\begin{lpic}[t(0mm),b(0mm),r(0mm),l(0mm)]{pics/key-lemma(0.7)}
\lbl[tr]{32,20;$w$}
\lbl[t]{39,18.5;$p_\delta$}
\lbl[t]{61,16;$p$}
\lbl[tr]{27,26;$q_\delta$}
\lbl[tr]{4,54;$q$}
\lbl[l]{35,24;$a_\delta$}
\lbl[lb]{41.5,31.5;$a$}
\lbl[br]{37.5,29.5,45;{\small $\alpha$}}
\lbl[W]{52,50;$\Omega_p$}
\lbl[W]{30,58;$\Omega_q$}
\end{lpic}
\end{wrapfigure}

Set $r_\delta=\dist{w}{a_\delta}{}$.
Consider the radial curve $\alpha\:[r_\delta,R]\to \Omega$ 
with respect to $w$ such that $\alpha(r_\delta)=a_\delta$.
Set $a=\alpha(R)$.
By Proposition~\ref{prop:rad-loc}, we have
\begin{align*}
\dist{p}{a}{}
&\le\dist{\~p}{\~a}{},
&
\dist{q}{a}{}
&\le\dist{\~q}{\~a}{}.
\end{align*}
Hence Key Lemma follows.\qeds

\begin{thm}{Lemma}\label{lem:good+good}
Let $\Omega_p$ and $\Omega_q$ be two $\kappa$-domains in $\bar X$.
Let 
$p\in X\cap \Omega_p$, 
$q\in X\cap\Omega_q$
and
$[pq]\subset X\cap(\Omega_p\cup\Omega_q)$.

Then for any geodesic $[qs]\subset \Omega_q\cap X$
the condition \ref{(**)} holds if $\dist{q}{\bar s}{}$ is sufficiently small. 
\end{thm}

\parit{Proof.}
Choose $w\in [pq]\cap\Omega_p\cap\Omega_q$.
Since $\Omega_q$ is a $\kappa$-domain, 
we have 
$$
\angle\hinge w{\bar s}q
\ge
\angk\kappa{w}{\bar s}{q}.
$$
for $\bar s\in \l]qs\r]$.
Therefore
$$\angle\hinge w{\bar s}p
\le
\pi-\angk\kappa{w}{\bar s}{q}.
$$
Note that for small values of $\dist{q}{\bar s}{}$,
we can apply Key Lemma;
hence the result.
\qeds

\begin{thm}{Corollary}\label{cor:good+good}
Let $\Omega_1$ and $\Omega_2$ be two $\kappa$-domains in $\bar X$.
Assume 
$$\Omega_3\z\subset \Omega_1\cup\Omega_2$$
is an open set such that for any two points $x,y\in X\cap \Omega_3$ any geodesic $[xy]$ lies in $\Omega_1\cup\Omega_2$.
Then $\Omega_3$ is a $\kappa$-domain.
\end{thm}

The following Lemma makes possible to produce triple of the domains 
$\Omega_1$, $\Omega_2$ and $\Omega_3$ as in the corollary above.

\begin{thm}{Lemma}\label{lem:xyz}
Let $[pq]$ be a geodesic in $X$
and the points $x$, $y$ and $z$ appear on $\l]pq\r[$
in the same order.
Assume that there are $\kappa$-domains
$\Omega_1\supset [xy]$ and $\Omega_2\supset [yz]$ in $\bar X$.
Then there is 
an open set $\Omega_3\subset \bar X$ which contains $[xz]$ and such that 
for any two points $v,w\in \Omega_3\cap X$ any geodesic $[vw]$ lies in $\Omega_1\cup\Omega_2$.

In particular, by Corollary~\ref{cor:good+good}, 
$\Omega_3$ is a $\kappa$-domain.
\end{thm}

Before the proof starts, 
let us discuss ultralimits of metric spaces briefly;
see \cite{akp} for more details.

Fix a nonprinciple ultrafilter $\omega$ on the natural numbers.
Denote by $\bar X^\omega$ the $\omega$-power of $\bar X$.

The space $\bar X$ will be considered as a subspace of $\bar X^\omega$ in the natural way.
Given a point $p\in \bar X$, we will denote by 
$B(p,\eps)_{\bar X}$ 
and $B(p,\eps)_{\bar X^\omega}$ 
 the $\eps$-ball centered at $p$ in $\bar X$ 
 and in  $\bar X^\omega$ correspondingly.

In is straightforward to check the following
\begin{clm}{}\label{clm:k-domain^o}
If $B(p,\eps)_{\bar X}$ is a $\kappa$-domain 
then so is $B(p,\eps)_{\bar X^\omega}$.
\end{clm}

\parit{Proof.} Arguing by contradiction,
assume that there is a sequence of geodesics $[u_nv_n]$ such that 
$u_n\to u\in [xz]$, $v_n\to v\in[xy]$ 
and $[u_nv_n]\not\subset\Omega_1\cup\Omega_2$ for each $n$.

The $\omega$-limit of $[u_nv_n]$
is a geodesic in $\bar X^\omega$ from $u$ to $v$ which does not lie in $[pq]$.
I.e., geodesics in $\bar X^\omega$ bifurcate at some point, say $w\in [xz]$.

According to \ref{clm:k-domain^o}, if $\eps>0$ is small enough, the ball $B(w,\eps)_{\bar X^\omega}$ forms a $\kappa$-domain, a contradiction. 
\qeds

\section{The proof}\label{sec:end}

\parit{Proof of \ref{({*})}.}
Note that one can split the geodesic $[pq]$
into segments 
in such a way that each segment lies in a $\kappa$-domain. 
More precisely, there is a sequence of points $p=p_0,p_1,\dots,p_n=q$
on $[pq]$ 
such that the sequence 
$$\dist{p}{p_0}{}, \dist{p}{p_1}{},\dots,\dist{p}{p_n}{}$$ 
is increasing and each geodesic 
$[p_{i-1}p_i]$ lies in a $\kappa$-domain.
Given $\eps>0$, the sequence above can be chosen in such a way that in addition
$\dist{p}{p_1}{}<\eps$.

Applying Lemma~\ref{lem:xyz} few times,
we get that the segment $[p_1p_{n-1}]$ belongs to one $\kappa$-domain.
Applying Key Lemma (\ref{lem:key}) to three points 
$p_1, p_{n-1}$ and $\bar s\in \l]qs\r]$ 
with small enough $\dist{q}{\bar s}{}$,
we get that \ref{({*})} holds for $p_1$, $[qs]$ and $\kappa_1=\kappa$.

Finally since $\dist{p}{p_1}{}$ can be made arbitrary small,
the triangle inequality implies that \ref{({*})} holds for $p$, $[qs]$ 
and arbitrary $\kappa_1<\kappa$.
\qeds

\section{Remarks}

\parbf{Cat's cradle construction.}
An alternative proof of the Key Lemma can be build on the so called \emph{Cat's cradle construction} 
from \cite{ab}.
This way you do not have to learn who are the radial curves and what is radial monotonicity.

\medskip

Choose small $\eps>0$ and apply the following procedure:
\begin{itemize}
\item Set $w_0=w$.
\item Choose $w_1\in [pw_0]$ so that $\dist{w_0}{w_1}{}=\eps$.
\item Choose $w_2\in [qw_1]$ so that $\dist{w_1}{w_2}{}=\eps$.
\item Choose $w_3\in [pw_2]$ so that $\dist{w_2}{w_3}{}=\eps$.
\item And so on.
\end{itemize}

\begin{wrapfigure}{r}{52mm}
\begin{lpic}[t(-3mm),b(3mm),r(0mm),l(0mm)]{pics/cradle(0.7)}
\lbl[rb]{1,37;$p$}
\lbl[lb]{70,37;$q$}
\lbl[lt]{35,0;$w_0=w$}
\lbl[b]{28,8.5;$w_1$}
\lbl[b]{39,14.5;$w_2$}
\lbl[b]{28,18.5;$w_3$}
\lbl[b]{39,23.5;$w_4$}
\lbl[b]{28,27.5;$w_5$}
\lbl[b]{40,33;$w_6$}
\end{lpic}
\end{wrapfigure}

Further, find a nice estimate for the value
$$\ell_n=\dist{p}{w_{2\cdot n}}{}+\dist{w_{2\cdot n}}{q}{}$$
in terms of 
$\angle\hinge wpq$,
$\dist{w}{p}{}$,
$\dist{w}{q}{}$
and
$$s_n\z=\sum_{i=1}^n\dist{w_{2\cdot(i-1)}}{w_{2\cdot i}}{}.$$
Note that 
$w_0$, $w_1,\z\dots,w_{2\cdot n}\in \Omega_p\cap\Omega_q$ 
if $s_n<R-\eps$.

If you were able to do everything as suggested,
you should get a weaker version of Key Lemma.
Namely you prove the required estimate if $R$ is bigger than bisector of $[\~p\~w\~q]$
at $\~w$.
This is good enough for the rest of the proof.
Playing a bit with the construction,
namely making $\eps$ depend on $n$,
you can actually get the Key Lemma in full generality.

(After making all this work you might appreciate the radial curves.)

\parbf{Finite dimensional case.}
In case $\bar X$ has finite (say Hausdorff) dimension,
one can build an easier proof using the formula of second variation see \cite{2nd-var}. 







\begin{thebibliography}{BFK}
\bibitem{ab} S. Alexander, R. Bishop, 
\textit{The Hadamard--Cartan theorem in locally convex spaces}, 
l'Enseignement Math., 36 (1990), 309--320.

\bibitem{akp}
S. Alexander, 
V. Kapovitch, 
A. Petrunin, 
\textit{Alexandrov geometry,} a draft available at \href{http://www.math.psu.edu/petrunin/}{\tt www.math.psu.edu/petrunin}

\bibitem{bgp}Burago, Yu.,  Gromov, M., and  Perelman, G.
\textit{A.D.~Aleksandrov
spaces with curvatures bounded below}, 
Uspehi Mat. Nauk 47: 2 (1992), 3-51, 222; translation
in Russian Math. Surveys 47: 2 (1992), 1-58.


\bibitem{2nd-var} Petrunin, A. \textit{Parallel transportation for Alexandrov space with curvature bounded below.}
Geom. Funct. Anal. 8 (1998), no. 1, 123--148.

\end{thebibliography}
\end{document}